\documentclass[11pt,a4paper]{amsart}
\usepackage[a4paper,hmargin=3cm,vmargin=3cm]{geometry}
\usepackage{amsmath,amsthm,amsfonts,amssymb}
\usepackage[latin1]{inputenc}
\usepackage[english]{babel}
\usepackage{xcolor}
\usepackage[hidelinks]{hyperref}
\AtBeginDocument{%
  \setlength{\abovedisplayskip}{0.8em plus .15em minus .15em}
  \setlength{\belowdisplayskip}{0.8em plus .15em minus .15em}
  \setlength{\abovedisplayshortskip}{0.6em plus .1em}
  \setlength{\belowdisplayshortskip}{0.6em plus .1em}
  \setlength{\parskip}{0.3em plus .1em minus .05em}
}
\makeatletter
\def\and{\unskip\space and\space}
\makeatother
\newcommand{\Nil}{\mathrm{Nil}_3}
\def\R{\mathbb{R}}
\def\H{\mathbb{H}}
\def\S{\mathbb{S}}
\def\cD{\mathcal{D}}
\def\cJ{\mathcal{J}}
\def\cA{\mathcal{A}}
\def\cV{\mathcal{V}}
\def\cF{\mathcal{F}}

\def\R{\mathbb{R}}

\newcommand{\ben}{\begin{enumerate}}
\newcommand{\bit}{\begin{itemize}}
\newcommand{\een}{\end{enumerate}}
\newcommand{\eit}{\end{itemize}}

\newcommand{\ed}{\end{document}}

\def\cA{\mathcal{A}}

\def\cD{\mathcal{D}}
\def\cR{\mathcal{R}}

\def\cV{\mathcal{V}}
\def\cJ{\mathcal{J}}

\def\cF{\mathcal{F}}

\let\8=\infty \let\0=\emptyset

\let\parc=\partial

\def\flecha{\rightarrow}
\def\esiz{\langle}
\def\esde{\rangle}

\def\cte.{\mathop{\rm cte.}\nolimits}

\def\det{\mathop{\rm det}\nolimits}

\def\R{\mathbb{R}}

\def\H{\mathbb{H}}
\def\S{\mathbb{S}}

\def\Nil{\rm{Nil}_3}

\newtheorem{mainthm}{Theorem}

\newtheorem{theorem}{Theorem}[section]
\newtheorem{lemma}[theorem]{Lemma}
\newtheorem{proposition}[theorem]{Proposition}
\newtheorem{remark}[theorem]{Remark}
\numberwithin{equation}{section}
\begin{document}
\title[The isoperimetric problem in Heisenberg $3$-space]{The isoperimetric problem in the Riemannian \\[0.1cm]Heisenberg $3$-space}


\author{Jos\'e M. Espinar}

\address{Departamento de Geometr\'{\i}­a y Topolog\'{\i}a, Universidad de Granada, Spain.}
\email{jespinar@ugr.es} 

\author{Pablo Mira}

\address{Departamento de Matem\'atica Aplicada y Estad\'{\i}stica, Universidad Polit\'ecnica 
de Cartagena, Spain} \email{pablo.mira@upct.es}

\thanks{Research supported by Project PID2024-160586NB-I00 funded by MCIN/AEI/10.13039/501100011033}




\begin{abstract}
We prove that every closed weakly stable constant mean curvature surface immersed in the Riemannian Heisenberg space $\Nil$ is an embedded rotational sphere. As a consequence, we solve the isoperimetric problem in $\Nil$: the isoperimetric domains of $\Nil$ are exactly the regions enclosed by rotational constant mean curvature spheres, and for each prescribed volume there is a unique such domain, up to ambient isometry.
\end{abstract}
\maketitle

\section{Introduction}
The isoperimetric problem in a Riemannian manifold $(M^n,g)$ asks to find, among the regions $\Omega\subset M$ of a given finite volume $V>0$, those that minimize the area of $\partial \Omega$. In this paper we solve this problem in the Riemannian Heisenberg space $\Nil$.

The study of the isoperimetric problem in $\Nil$ goes back to Tomter in 1993. In \cite{To}, he proved that for any $V>0$ there exists a constant mean curvature (CMC) rotational sphere embedded in $\Nil$, enclosing a region of volume $V$, and studied their isoperimetric properties. In addition, Tomter proved \cite{To} that, for small volumes, these rotational spheres are isoperimetric; see also Nardulli \cite{Na1,Na2}. A known conjecture in the theory asserts then that all isoperimetric regions of $\Nil$ are bounded by rotational spheres; see e.g. \cite{BT,DHM,FM,FMP,MMPR14,MP,MPT}. In this work we give a positive answer to this conjecture for all volumes.

\begin{mainthm}\label{thm:iso}
The isoperimetric domains of $\Nil$ are exactly the bounded regions enclosed by rotational constant mean curvature spheres. For each $V>0$, the isoperimetric domain of volume $V$ is unique up to ambient isometry.
\end{mainthm}

Theorem \ref{thm:iso} will follow from the classification of closed weakly stable CMC surfaces in $\Nil$ in Theorem \ref{thm:main} below. To state this result, recall that constant mean curvature surfaces in a Riemannian $3$-manifold are critical points of the area functional under volume-preserving variations. For a closed two-sided CMC surface $\Sigma$, the second variation is governed by the Jacobi operator
\begin{equation}\label{jacop}
L=\Delta+|S|^2+\operatorname{Ric}(N),
\end{equation}
where $\Delta,S,N$ denote, respectively, the Laplacian, the shape operator and the unit normal of $\Sigma$. The associated index form is
\begin{equation}\label{eq:Q}
Q(v,w)=\int_\Sigma\bigl(\langle\nabla v,\nabla w\rangle
-(|S|^2+\operatorname{Ric}(N))\, vw\bigr)\,dA=-\int_\Sigma vLw\,dA.
\end{equation}
We say that $\Sigma$ is \emph{weakly stable} if
\begin{equation}\label{wike}
Q(v,v)\geq0\qquad\text{for every }v\in C^\infty(\Sigma)
\text{ such that }\int_\Sigma v\,dA=0.
\end{equation}
Barbosa and do Carmo \cite{BC} proved that the only closed weakly stable CMC hypersurfaces in Euclidean space $\R^{n+1}$ are round spheres. Also, in \cite{BCE} Barbosa, do Carmo and Eschenburg established the corresponding classification in $\S^{n+1}$ and $\H^{n+1}$. In $\Nil$, Torralbo and Urbano \cite{TUstable} proved that all rotational CMC spheres are weakly stable, and that a closed weakly stable CMC surface with $H^2\geq3$ is either a sphere or an embedded torus. For further results on weak stability in homogeneous $3$-manifolds, see Nelli--Rosenberg \cite{NR}, Souam \cite{So}, and Torralbo--Urbano \cite{TUstable}.

Our main result extends the Barbosa--do Carmo theorem to $\Nil$:

\begin{mainthm}\label{thm:main}
Any closed weakly stable CMC surface immersed in $\Nil$ is an embedded rotational CMC sphere.
\end{mainthm}
Let us explain how Theorem \ref{thm:iso} follows from Theorem \ref{thm:main}. It is well-known that if $\Omega$ is an isoperimetric region of $\Nil$, then $\partial \Omega$ is a closed, weakly stable CMC surface. Since $\Nil$ is a homogeneous manifold diffeomorphic to $\R^3$, a classical argument proves that $\parc\Omega$ is connected, see e.g. \cite{DanielSurvey,MMPR14}. Thus, Theorem \ref{thm:main} indeed implies that any isoperimetric region in $\Nil$ is bounded by some rotational CMC sphere. On the other hand, since $\Nil$ is homogeneous, isoperimetric regions exist for every prescribed volume $V>0$; see e.g. \cite{Ros}. Hence, the remaining statements in Theorem \ref{thm:iso} follow then from Tomter \cite{To}, who proved that for any $V>0$ there is exactly one rotational CMC sphere enclosing that volume, up to ambient isometry.

\subsection*{Strategy of the proof}
A well-known theorem by Abresch-Rosenberg \cite{AR2,AR05} proves that any CMC sphere in a rotationally symmetric homogeneous $3$-manifold, and in particular in $\Nil$, is a sphere of revolution; see also \cite{MMPR1,MMPR2} for the classification of CMC spheres in arbitrary homogeneous $3$-manifolds. Thus, the idea of the proof is to show that there are no closed weakly stable CMC surfaces in $\Nil$ of genus $\gamma\geq 1$.

Arguing by contradiction, let $\Sigma$ be any such surface, and denote for any $f\in C^{\8}(\Sigma)$ $$
[f]:=\frac{1}{|\Sigma|}\int_\Sigma f\,dA, \hspace{0.5cm} |\Sigma|:=\operatorname{Area}(\Sigma).
$$
Consider the Heisenberg dilation field $\cD=x\partial_x+y\partial_y+2z\partial_z$, and denote $\phi:=\langle\cD,N\rangle$, where $N$ is the unit normal of $\Sigma$. Then we will use $\phi_0:=\phi-[\phi]$ as a test function for the stability inequality \eqref{wike}. In $\R^3$, Barbosa-do Carmo proved that the analogous function $\phi_0$ satisfies $Q(\phi_0,\phi_0)<0$ unless $\Sigma$ is totally umbilic. This is not true in $\Nil$. We will use instead an alternative argument based on integral formulas and a Cauchy-Schwarz inequality to prove that the condition $Q(\phi_0,\phi_0)\geq 0$ implies $[\nu^2]<1/2$, where $\nu:=\esiz N,\parc_z\esde$ is the usual vertical angle function for surfaces in $\Nil$.

Once here, we use a second test function on \eqref{wike}, namely, 
\begin{equation}\label{fh}
\begin{gathered}
v_H:=f_H(\nu)-\frac{[f_H(\nu)]}{[\phi]}\phi, \hspace{0.5cm} f_H(t):=1+\frac{t}{H}\arctan\frac{t}{H}.
\end{gathered}
\end{equation}
The function $f_H(\nu)$ in \eqref{fh} was used by Torralbo-Urbano \cite[Section 4]{TUstable} on rotational CMC spheres of $\Nil$ in order to prove their weak stability. We will use it instead applied on any weakly stable CMC surface of genus $\gamma\geq 1$, and obtain using the previous inequality $[\nu^2]<1/2$ that
\begin{equation}\label{inegen}
Q(v_H,v_H)< 6\pi (1-\gamma)\leq 0,
\end{equation} 
a contradiction.

The paper is organized as follows. In Section \ref{sec:prelim} we provide some preliminary equations and prove a Minkowski-type formula for closed oriented surfaces in $\Nil$. In Section \ref{sec:cmc} we give some integral formulas for closed CMC surfaces in $\Nil$, using divergence identities and first and second variation formulas. In Section \ref{sec:proofB} we prove Theorem \ref{thm:main}, in three steps. First we use $\phi_0$ to obtain the estimate $[\nu^2]<1/2$ on any closed weakly stable CMC surface. In Step 2 we compute $Q(v_H,v_H)$ in terms of $[f_H(\nu)]$. In Step 3, we finally assume that $\Sigma$ has genus $\gamma\geq 1$, and use the previous steps to arrive at  the inequality \eqref{inegen}, contradicting weak stability. This implies $\gamma=0$, and the Abresch--Rosenberg theorem \cite{AR05} completes the proof. 

\section{A Minkowski-type formula}\label{sec:prelim}

Up to homothety, $\Nil$ can be seen as $\R^3$ endowed with the Riemannian metric $ds^2=dx^2+dy^2+\bigl(dz+y\,dx-x\,dy\bigr)^2$. This metric is left-invariant for the usual group structure of $\Nil\equiv \H^1$, given by 
$$
(x,y,z)\cdot (x',y',z') =\bigl(x+x',y+y',z+z'+ xy'-yx'\bigr).$$ Its isometry group is $4$-dimensional and admits $\S^1$-rotational subgroups. We work in $\Nil=\mathbb E(0,1)$ in Daniel's notation \cite{Daniel}. In coordinates $(x,y,z)$ the left-invariant orthonormal frame is
\begin{equation}\label{frame1}
E_1=\partial_x-y\partial_z,\qquad
E_2=\partial_y+x\partial_z,\qquad E_3=\partial_z=\xi.
\end{equation}
We orient $\Nil$ so that $\{E_1,E_2,E_3\}$ is positively oriented.

Let $\Sigma$ be a closed immersed oriented surface in $\Nil$, with unit normal $N$ and shape operator $SY=-\overline\nabla_YN$. We write $H=\tfrac12\operatorname{tr}S$ for its mean curvature, and set $JY=N\times Y$. Let $
\nu:=\langle N,\xi\rangle$ be the \emph{angle function} of $\Sigma$, denote $T:=\xi-\nu N$, and observe that $|T|^2+\nu^2=1.$ 
The compatibility equations in \cite[Proposition 3.3]{Daniel} give
\begin{equation}\label{eq:Daniel}
\nabla_YT=\nu(SY-JY),\qquad
Y(\nu)+\langle SY-JY,T\rangle=0.
\end{equation}
In particular,
\begin{equation}\label{eq:basic}
\nabla\nu=-ST-JT,\qquad
\operatorname{div}T=2H\nu,\qquad
\operatorname{div}(JT)=2\nu.
\end{equation}
For the last identity, take the trace of
$\nabla_Y(JT)=\nu(JSY+Y)$ and use $\operatorname{tr}(JS)=0$.
The Gauss equation is (see \cite[Corollary 3.2]{Daniel})
\begin{equation}\label{eq:Gauss}
K=\det S+1-4\nu^2.
\end{equation}

Consider next the Heisenberg dilations
\begin{equation}\label{dila}
\delta_t(x,y,z)=(e^tx,e^ty,e^{2t}z)
\end{equation}
and their velocity field $\cD=x\partial_x+y\partial_y+2z\partial_z$.
We set  $\Sigma_t:=\delta_t(\Sigma)$, and write $\cA(t)$ and $\cV(t)$ for its area and algebraic volume. The volume form in $\Nil$ is $dx\,dy\,dz$, and hence $\cV(t)=e^{4t}\cV (0)$, where $\cV(0):=V$ is the algebraic volume of $\Sigma$. The area element of $\Sigma_t$ varies according to $\delta_t^* (dA_t)={\rm Jac} (\delta_t) \, dA$. In addition, the area Jacobian of $\delta_t$ is (see \cite[Proof of Thm. 4.12]{RR}) $$\cJ(t):={\rm Jac} (\delta_t)=e^{2t}\sqrt{\nu^2 + e^{2t} (1-\nu^2)}.$$ Thus, $\cJ'(0)=3-\nu^2$ and $\cJ''(0)=9-4\nu^2-\nu^4$, and so 
\begin{equation}\label{dosa}
\cA'(0)=\int_\Sigma (3-\nu^2)\, dA, \hspace{0.5cm} \cA''(0)= \int_\Sigma  (9-4\nu^2-\nu^4)\, dA.
\end{equation}
For the volume,
\begin{equation}\label{dosv}
\cV'(0)=4 \cV(0), \hspace{0.5cm} \cV''(0)=16 \cV(0).
\end{equation}
As a consequence we have:
\begin{lemma}[Minkowski-type formula]\label{lem:minkowski}
For every closed immersed oriented surface $\Sigma$ in $\Nil$, denoting $\phi:=\langle \cD,N\rangle$, we have
\begin{equation}\label{integ1}
2\int_\Sigma H\phi\,dA=-\int_\Sigma(3-\nu^2)\,dA.
\end{equation}
\end{lemma}
\begin{proof}
The variation $\Sigma_t=\delta_t(\Sigma)$ has initial normal velocity $\phi$. Thus, we have the well-known formulas 
\begin{equation}\label{avo}
\cA'(0)=-2\int_\Sigma H \phi \, dA, \hspace{0.5cm} \cV'(0)=-\int_\Sigma \phi \, dA.
\end{equation}
Here, we have chosen the sign of the algebraic volume so that it is positive for embedded surfaces with their inner orientation. 
From \eqref{dosa} and the first equation in \eqref{avo}, we obtain directly \eqref{integ1}.
\end{proof}

\section{Integral formulas for CMC surfaces in \texorpdfstring{$\Nil$}{Nil3}}\label{sec:cmc}
Let us assume next that $\Sigma$ has constant mean curvature $H$. Since there are no closed minimal surfaces in $\Nil$, we suppose up to a change of orientation that $H>0$. Thus, $\Sigma$ is a critical point of the functional ${\rm Area} -2H\, {\rm Vol}$. The Hessian of this functional is the index form $Q$ in \eqref{eq:Q}, evaluated on the normal components of the variation, see e.g. \cite[Proposition 2.5]{BCE}.
The next lemma computes $Q$ on the dilation function $\phi=\langle \mathcal{D},N\rangle$ and on constants. We will use these identities to test weak stability in Section \ref{sec:proofB}.

\begin{lemma}\label{lem:dilation-index}
Let $\Sigma$ be a closed immersed oriented CMC surface in $\Nil$. Then
\begin{align}
Q(\phi,\phi)&=-\int_\Sigma(3+\nu^4)\,dA,\label{sebau3}\\
Q(1,1)&=-4H^2|\Sigma|+4\int_\Sigma\nu^2\,dA
+4\pi\chi(\Sigma).\label{q11}
\end{align}
Also, for every $\varphi\in C^\infty(\Sigma)$, we have
\begin{equation}\label{quv}
Q(\phi,\varphi)=2\int_\Sigma
\bigl(\nu\langle\nabla \varphi,T\rangle+H \varphi(1+\nu^2)\bigr)\,dA.
\end{equation}
\end{lemma}
\begin{proof}
Let us first note that, by \eqref{dosv}, \eqref{avo}, the algebraic volume of $\Sigma$ is $V=-\tfrac14\int_\Sigma\phi\,dA$, and then, by \eqref{integ1},
\begin{equation}\label{voli}
[\phi]=-\frac{3-[\nu^2]}{2H}<0,\qquad
8HV=\int_\Sigma(3-\nu^2)\,dA.
\end{equation}
We now apply the second variation formula to the variation $\Sigma_t=\delta_t(\Sigma)$, to obtain using \eqref{dosa}, \eqref{dosv} and \eqref{voli} that
$$
\begin{aligned}
Q(\phi,\phi)&=\cA''(0)-2H\cV''(0)=\int_\Sigma(9-4\nu^2-\nu^4)\,dA-32HV\\
&=-\int_\Sigma(3+\nu^4)\,dA,
\end{aligned}$$
what proves \eqref{sebau3}. 

In order to prove \eqref{q11}, we first note that the Gauss equation \eqref{eq:Gauss}, together with $\operatorname{Ric}(N)=-2+4\nu^2$ and $|S|^2=4H^2-2\det S$, show that the Jacobi operator \eqref{jacop} of $\Sigma$ is given by
\begin{equation}\label{estapot}
L=\Delta +P, \hspace{0.5cm} P:=4H^2-2K-4\nu^2.
\end{equation}
Then, \eqref{q11} follows directly from \eqref{eq:Q}, \eqref{estapot} and Gauss--Bonnet.

Finally, to prove \eqref{quv}, take $\varphi\in C^{\infty}(\Sigma)$ and consider the two-parameter variation of $\Sigma$ given by $$\Phi_{s,t}(p)=\delta_t\bigl(\exp_p(s\,\varphi(p)N(p))\bigr):\Sigma\flecha \Nil.$$ Here, $\delta_t$ denotes the Heisenberg dilation in \eqref{dila} and $\exp$ is the exponential map in $\Nil$. Consider the functional
\[
\cF(s,t):=\operatorname{Area}(\Phi_{s,t}(\Sigma))
-2H\operatorname{Vol}(\Phi_{s,t}(\Sigma)),
\]
with $H$ fixed. By \cite[Prop. 2.5]{BCE}, $Q(\phi,\varphi)=\partial_{st} \cF(0,0)$.  Applying the dilation formulas \eqref{dosa}--\eqref{dosv} to $\Sigma_s:=\Phi_{s,0}(\Sigma)$, we obtain
\begin{equation}\label{fs}
\partial_t\cF(s,0)
=
\int_{\Sigma}(3-\nu_s^2)\,dA_s
-8H\operatorname{Vol}(\Sigma_s),
\end{equation}
where $dA_s$ and $\nu_s$ denote the area measure and the angle function of $\Sigma_s$, respectively, pulled back to $\Sigma$ by $\Phi_{s,0}$.

We next differentiate \eqref{fs} with respect to $s$. A dot denotes differentiation (covariant differentiation for vector fields) with respect to $s$ at $(s,t)=(0,0)$. Since $\partial_s \Phi_{s,0}=\varphi N$ at $s=0$, the standard normal variation formulas give
\begin{equation}\label{dot1}
\dot{\left(dA_s\right)}= -2H \varphi \,dA, \hspace{0.5cm} \dot{\left ({\rm Vol}(\Sigma_s)\right)}=-\int_\Sigma \varphi \, dA.
\end{equation}
Also, using $\dot\Phi = \varphi N$, 
\begin{equation}\label{dot2}
\dot\nu = \langle \dot{N},\xi \rangle + \langle N, \overline\nabla_{\varphi N} \xi\rangle = -\langle \nabla \varphi,\xi\rangle + \varphi\langle N, \overline\nabla_N \xi\rangle = -\langle \nabla \varphi,T\rangle,
\end{equation} 
where in the last step we used that $\xi$ is Killing and its tangential part is $T$. Differentiating \eqref{fs} and using \eqref{dot1}--\eqref{dot2}, we obtain
\begin{equation}\label{dot3}
\begin{aligned}
\dot{\big(\partial_t \cF(s,0)\big)}
&=\int_\Sigma
\left(2\nu\langle\nabla \varphi,T\rangle
-2H \varphi(3-\nu^2)+8H \varphi\right)dA\\
&=2\int_\Sigma
\left(\nu\langle\nabla \varphi,T\rangle+H \varphi(1+\nu^2)\right)dA.
\end{aligned}
\end{equation}
Thus, since $Q(\phi,\varphi)= \partial_{st} \cF(0,0)=\dot{\big(\partial_t \cF(s,0)\big)}$, \eqref{dot3} gives \eqref{quv}.
\end{proof}

We finish this section with an integral formula for the angle function. 
\begin{lemma}\label{lem:even}
Let $\Sigma$ be a closed oriented CMC surface immersed in $\Nil$, let $f\in C^\infty[-1,1]$ be even, and set $F(t)=\int_0^t f(s)\,ds$. Then
\begin{equation}\label{foreven}
\int_\Sigma(K+4\nu^2)f(\nu)\,dA
=4\int_\Sigma\nu F(\nu)\,dA
+2\pi\chi(\Sigma)F(1).
\end{equation}
\end{lemma}
\begin{proof}
The divergence formula in Torralbo--Urbano \cite[proof of Theorem 1]{TUgauss}, specialized to $(\kappa,\tau)=(0,1)$, yields
\begin{equation}\label{eq:TU}
{\rm div} (\nabla_TT-2H\nu T)=K(1-3\nu^2)-8\nu^4.
\end{equation}
Also, applying \eqref{eq:Daniel} to $X=T$ gives
\begin{equation}\label{ecuba1}
\nabla_T T -2 H\nu T = -\nu W, \hspace{0.5cm} W:=(2H I-S+J)T.
\end{equation}
Putting together \eqref{eq:TU} and \eqref{ecuba1} we obtain
\begin{equation}\label{eq:divnuW}
\operatorname{div}(\nu W)=-K(1-3\nu^2)+8\nu^4.
\end{equation}
Let $a\in C^\infty[-1,1]$. To compute
$\operatorname{div}\bigl(\nu a(\nu)W\bigr)$, we first evaluate,
using \eqref{eq:basic} and the identities
$S^2-2HS+(\det S)I=0$ and $SJ+JS=2HJ$,
the gradient term
\begin{equation}\label{eq:contraction}
\langle\nabla\nu,W\rangle =-\langle T,(S-J)(2HI-S+J)T\rangle =-(K+4\nu^2)(1-\nu^2);
\end{equation}
here the last equality also uses the Gauss equation \eqref{eq:Gauss}. Then, the product rule and \eqref{eq:divnuW}-\eqref{eq:contraction} yield
\[
\operatorname{div}\bigl(\nu a(\nu)W\bigr)
=4\nu G(\nu)-(K+4\nu^2)G'(\nu), \hspace{0.5cm} G(t):=t(1-t^2)a(t).
\]
Thus, integrating, we obtain
\begin{equation}\label{intege}
\int_\Sigma(K+4\nu^2)G'(\nu)\,dA
=4\int_\Sigma\nu G(\nu)\,dA,
\end{equation}
for any function $G(t)=t(1-t^2)a(t)$, with $a\in C^\infty[-1,1]$.

We want $f,F$ to appear in \eqref{intege}. For that, we choose $G(t):=F(t)-tF(1)$. Since $f$ is even, $F(-1)=-F(1)$ and so $G$ vanishes at $0,\pm 1$. Hence, $a(t):=G(t)/(t(1-t^2))$ extends smoothly to $[-1,1]$. Thus, we can apply \eqref{intege} with our $G(t)$, giving
\[
\int_\Sigma(K+4\nu^2)\bigl(f(\nu)-F(1)\bigr)\,dA
=4\int_\Sigma\nu\bigl(F(\nu)-\nu F(1)\bigr)\,dA.
\]
The terms involving $4\nu^2F(1)$ cancel, and Gauss--Bonnet proves \eqref{foreven}.
\end{proof}
\begin{remark}
The potential of the Jacobi operator $L=\Delta +P$ in \eqref{estapot} can be rewritten as
$
P=4(H^2+\nu^2)-2(K+4\nu^2)$. Thus, \eqref{foreven} can be written alternatively as
\begin{equation}\label{forevenP}
\int_\Sigma P f(\nu)\,dA
=
4\int_\Sigma
\bigl((H^2+\nu^2)f(\nu)-2\nu F(\nu)\bigr)\,dA
-4\pi\chi(\Sigma)F(1),
\end{equation}
with $f\in C^\infty[-1,1]$ even, and $F(t)=\int_0^t f(s)\,ds$.
\end{remark}

\section{Proof of Theorem B}\label{sec:proofB}
Let $\Sigma$ be a closed, weakly stable CMC surface in $\Nil$, oriented so that $H>0$. We assume the genus of $\Sigma$ is $\gamma\geq1$, and so $\chi(\Sigma)=2-2\gamma\leq0$.

\subsection*{Step 1: The first stability test}
In this step we apply the weak stability condition $Q(\phi_0,\phi_0)\geq 0$ to the test function $\phi_0:=\phi-[\phi]$, to produce some intermediate bounds needed later on.
\begin{proposition}\label{pro:first}
In the above conditions, we have
\begin{equation}\label{coarsebounds}
[\nu^2]<\frac12.
\end{equation}
\end{proposition}
\begin{proof}
Since $[\phi_0]=0$, weak stability, \eqref{voli} and Lemma \ref{lem:dilation-index} give
\begin{equation}\label{firsttest}
0\leq\frac{H^2Q(\phi_0,\phi_0)}{|\Sigma|}
=(3-[\nu^2])^2\left([\nu^2]+\frac{\pi\chi(\Sigma)}{|\Sigma|}\right)
-H^2\bigl(6-10[\nu^2]+3[\nu^2]^2+[\nu^4]\bigr).
\end{equation}
Let us also note the elementary inequalities
\begin{equation}\label{moments}
[\nu^2]^2\leq[\nu^4]\leq[\nu^2]<1.
\end{equation}
The first two follow from Cauchy--Schwarz and $|\nu|\leq1$. The last one is a direct consequence of the fact that the horizontal distribution in $\Nil$ spanned by $\{E_1,E_2\}$ is not integrable. Thus, by \eqref{moments}, the last term in \eqref{firsttest} can be bounded as 
\begin{equation}\label{mom2}
\begin{aligned}
6-10[\nu^2]+3[\nu^2]^2+[\nu^4]
& =2(3-2[\nu^2])(1-[\nu^2])+[\nu^4]-[\nu^2]^2 \\ & \geq 2\big(3-2[\nu^2])(1-[\nu^2]\big).
\end{aligned}
\end{equation}
Plugging \eqref{mom2} into \eqref{firsttest} we obtain 
\begin{equation}\label{firsttestbound}
0<2H^2(3-2[\nu^2])(1-[\nu^2])
\leq(3-[\nu^2])^2\left([\nu^2]+\frac{\pi\chi(\Sigma)}{|\Sigma|}\right).
\end{equation}
In particular, 
\begin{equation}\label{anglelowerbound}
[\nu^2]+\frac{\pi\chi(\Sigma)}{|\Sigma|}>0.
\end{equation}
We next turn to prove that $[\nu^2]<1/2$. First, Lemma \ref{lem:even} applied to $f(t)=1-t^2$ gives
$$
[(K+4\nu^2)(1-\nu^2)]
=
\frac43\left(
3[\nu^2]-[\nu^4]
+\frac{\pi\chi(\Sigma)}{|\Sigma|}
\right).$$
Consequently, \eqref{anglelowerbound} and \eqref{moments} imply
\begin{equation}\label{weightedcurvature}
\bigl[(K+4\nu^2)(1-\nu^2)\bigr]>\frac43\bigl(2[\nu^2]-[\nu^4]\bigr)
\geq\frac43[\nu^2]>0.
\end{equation}
This positivity allows us to estimate the norm of $(S-HI)T$. Indeed, since $S-HI$ is symmetric and traceless, the Gauss equation \eqref{eq:Gauss} gives
\[
(S-HI)^2=(H^2-\det S)I=(H^2+1-K-4\nu^2)I.
\]
Thus, by \eqref{weightedcurvature},
\begin{equation}\label{cs1}
\begin{aligned}
\bigl[|(S-HI)T|^2\bigr]
&=(H^2+1)(1-[\nu^2])-[(K+4\nu^2)(1-\nu^2)]\\
&<(H^2+1)(1-[\nu^2]),
\end{aligned}
\end{equation}
Moreover, since \(T\perp JT\), we have
\begin{equation}\label{cs2}
|HT+JT|^2=(H^2+1)(1-\nu^2).
\end{equation}
Now, \eqref{eq:basic} and \eqref{cs2} yield
\[
\begin{aligned}
\operatorname{div}\bigl(\nu(HT+JT)\bigr)
&=\langle\nabla\nu,HT+JT\rangle
  +\nu\bigl(H\operatorname{div}T+\operatorname{div}(JT)\bigr)\\
&=-\langle(S-HI)T,HT+JT\rangle
  -|HT+JT|^2+2(H^2+1)\nu^2\\
&=-\langle(S-HI)T,HT+JT\rangle
  +(H^2+1)(3\nu^2-1).
\end{aligned}
\]
Integrating, we obtain
\[
[\langle(S-HI)T,HT+JT\rangle]=(H^2+1)(3[\nu^2]-1).
\]
From here, Cauchy--Schwarz and \eqref{cs1}-\eqref{cs2} give now
\[
\begin{aligned}
(H^2+1)^2(3[\nu^2]-1)^2
&\leq[|(S-HI)T|^2]\,[|HT+JT|^2]\\
&<(H^2+1)^2(1-[\nu^2])^2.
\end{aligned}
\]
Therefore, $[\nu^2]<1/2$. 
\end{proof}

\subsection*{Step 2: The test function $v_H$}
We now compute the index form of the function $v_H$ defined in \eqref{fh}. Since $Q(\phi,\phi)$ was already computed in \eqref{sebau3}, we merely need to find $Q(\phi,f_H(\nu))$ and $Q(f_H(\nu),f_H(\nu))$. We do this next.

\begin{lemma}\label{lem:energies}
The following identities hold:
\begin{equation}\label{entries2}
\begin{aligned}
Q(\phi,f_H(\nu))
&=2|\Sigma|\, H \bigl(\left(H^2+1\right) [f_H(\nu)]-H^2\bigr),\\[0.2cm]
Q(f_H(\nu),f_H(\nu))
&=-4|\Sigma| \,H^2\, [f_H(\nu)]-4\pi\chi(\Sigma)\,C_H,
\end{aligned}
\end{equation}
where 
\begin{equation}\label{CHval}
C_H:=-\frac34-\frac{H^2+3}{4H}\arctan\frac1H<0.
\end{equation}
\end{lemma}
\begin{proof}
We begin with the first equation in \eqref{entries2}. Applying \eqref{quv} to $\varphi=f_H(\nu)$ gives
\begin{equation}\label{qe1}
Q(\phi,f_H(\nu))
=2\int_\Sigma\bigl(\nu\langle\nabla f_H(\nu),T\rangle
+H(1+\nu^2)f_H(\nu)\bigr)\,dA.
\end{equation} 
Now, recalling the definition of $f_H$ in \eqref{fh} and denoting $F_H(t):=\int_0^t f_H(s)\, ds$, one has directly $(tf_H(t)-F_H(t))'=tf_H'(t)$. Since $\operatorname{div}T=2H\nu$ by \eqref{eq:basic}, integration by parts yields 
\begin{equation}\label{qe2}
\begin{aligned}
\int_\Sigma\nu\langle\nabla f_H(\nu),T\rangle\,dA
&=\int_\Sigma\langle\nabla(\nu f_H(\nu)-F_H(\nu)),T\rangle\,dA\\
&=-2H\int_\Sigma\bigl(\nu^2f_H(\nu)-\nu F_H(\nu)\bigr)\,dA.
\end{aligned}
\end{equation}
Therefore, from \eqref{qe1}-\eqref{qe2},
\begin{equation}\label{qe3}
Q(\phi,f_H(\nu))
=2H\int_\Sigma\bigl((1-\nu^2)f_H(\nu)+2\nu F_H(\nu)\bigr)\,dA.
\end{equation}
Now, using the definition of $f_H(t)$ in \eqref{fh}, a direct computation gives 
\begin{equation}\label{mixed}
(1-t^2)f_H(t)+2tF_H(t)=(H^2+1)f_H(t)-H^2.
\end{equation}
Plugging \eqref{mixed} into \eqref{qe3} proves the first identity in \eqref{entries2}.

We next prove the second identity in \eqref{entries2}. Write the Jacobi operator $L=\Delta + P$ as in \eqref{estapot}. Then, by \eqref{eq:Q},
\begin{equation}\label{q21}
Q(f_H(\nu),f_H(\nu))
=\int_\Sigma
\left(f_H'(\nu)^2|\nabla\nu|^2-Pf_H(\nu)^2\right)dA.
\end{equation}
Now, $\nu=\esiz N,\xi\esde$ is a Jacobi function, since $\xi$ is a Killing field of $\Nil$. Thus, $\Delta\nu=-P\nu$. Integration by parts in the gradient term of \eqref{q21} gives then
\begin{equation}\label{q22}
\begin{aligned}
Q(f_H(\nu),f_H(\nu))
&=-\int_\Sigma
\left(\int_0^\nu f_H'(s)^2\,ds\right)\Delta\nu\,dA
-\int_\Sigma Pf_H(\nu)^2\,dA\\
&=\int_\Sigma P\left(
\nu\int_0^\nu f_H'(s)^2\,ds-f_H(\nu)^2
\right)dA.
\end{aligned}
\end{equation}
Let us evaluate the inner expression. By definition of \(f_H\) in \eqref{fh}, a computation gives
$$
t\int_0^t f_H'(s)^2\,ds-f_H(t)^2
=-\frac32f_H(t)+\frac{H^2}{2(H^2+t^2)}.$$
Hence, by \eqref{q22},
\begin{equation}\label{q23}
Q(f_H(\nu),f_H(\nu))
=\int_\Sigma P\,R_H(\nu)\,dA, \hspace{0.5cm} R_H(t):=-\frac32f_H(t)+\frac{H^2}{2(H^2+t^2)}.
\end{equation}
To evaluate the integral in \eqref{q23}, we apply \eqref{forevenP} for $f(t):=R_H(t)$, which is a smooth even function. First, note that
\[
\mathcal R_H(t)
:=\int_0^tR_H(s)\,ds
=-\frac{3t}{4}
-\frac{H^2+3t^2}{4H}\arctan\frac tH,
\]
and in particular $\mathcal R_H(1)=C_H$, with $C_H$ as in \eqref{CHval}. So, by \eqref{q23}, \eqref{forevenP}, 
\begin{equation}\label{qe25}
Q(f_H(\nu),f_H(\nu))=4\int_\Sigma
\left((H^2+\nu^2)R_H(\nu)
-2\nu\mathcal R_H(\nu)\right)dA -4\pi\chi(\Sigma)C_H.
\end{equation}
Finally, direct substitution of the expressions of $f_H, R_H, \cR_H$ gives
\begin{equation}\label{qe24}
(H^2+t^2)R_H(t)-2t\mathcal R_H(t)
=-H^2f_H(t).
\end{equation}
Inserting \eqref{qe24} into \eqref{qe25} we obtain the second identity in \eqref{entries2}.
\end{proof}

Once here, we are in the position to compute $Q(v_H,v_H)$ for the function $v_H$ in \eqref{fh}. First, by \eqref{voli}, and denoting for simplicity $M:=[f_H(\nu)]$, we can write $v_H$ as \[
v_H=f_H(\nu)-\frac{M}{[\phi]}\phi
=f_H(\nu)+\frac{2HM}{3-[\nu^2]}\phi.
\]
We note that $[v_H]=0$ and $M\geq 1$, since $f_H(t)\geq 1$. 
Then, from \eqref{sebau3} and \eqref{entries2}, we have by bilinearity
\[
\begin{aligned}
\frac{Q(v_H,v_H)}{4H^2|\Sigma|}
={}&-M-\frac{\pi\chi(\Sigma)}{H^2|\Sigma|}C_H\\
&+\frac{2M}{3-[\nu^2]}
\bigl((H^2+1)M-H^2\bigr)
-\frac{3+[\nu^4]}{(3-[\nu^2])^2}M^2.
\end{aligned}
\]
We are interested in estimating the above quantity. For that, we first rearrange it by writing $(H^2+1)M-H^2=M+H^2(M-1)$, obtaining
\begin{equation}\label{qvv}
\begin{aligned}
\frac{Q(v_H,v_H)}{4H^2|\Sigma|}
={}&-M+\frac{3-2[\nu^2]-[\nu^4]}{(3-[\nu^2])^2}M^2\\
&+\frac{2H^2}{3-[\nu^2]}M(M-1)
-\frac{\pi\chi(\Sigma)}{H^2|\Sigma|}C_H.
\end{aligned}
\end{equation}
Since $M(M-1)\geq 0$ and $[\nu^2]<1/2$, we obtain from \eqref{qvv} the estimate 
\begin{equation}\label{roughenergy}
\frac{Q(v_H,v_H)}{4H^2|\Sigma|}
\leq p_H(M)-\frac{\pi\chi(\Sigma)}{H^2|\Sigma|}C_H,
\end{equation}
where
\begin{equation}\label{polipe}
p_H(x):=-x+\frac{x^2}{3}+\frac{4H^2}{5}x(x-1).
\end{equation}

\subsection*{Step 3: The final estimate}
Our final objective is to prove that the right-hand side in \eqref{roughenergy} is negative. For that, we will need the following estimate.
\begin{lemma}\label{lem:meanbound}
The mean value $M=[f_H(\nu)]$ satisfies
\begin{equation}\label{meanb2}
(H^2+1)(M-1)^2-(M-1)
\leq
-\frac{\pi\chi(\Sigma)}{2H^3|\Sigma|}
\arctan\frac1H.
\end{equation}
\end{lemma}
\begin{proof}
Set
\[
G_H(t):=\frac1H\arctan\frac tH.
\]
Since $G_H'(t)=(H^2+t^2)^{-1}$ and
$tG_H(t)=f_H(t)-1$, we can compute using \eqref{eq:basic}
\[
\operatorname{div}\bigl(G_H(\nu)(HT+JT)\bigr)
=
\frac{\langle\nabla\nu,HT+JT\rangle}{H^2+\nu^2}
+2(H^2+1)(f_H(\nu)-1).
\]
Integrating, we have
\begin{equation}\label{cs0}
\left[\frac{\langle\nabla\nu,HT+JT\rangle}{H^2+\nu^2}\right]
=-2(H^2+1)(M-1).
\end{equation}
From \eqref{cs2}, we can apply Cauchy--Schwarz in \eqref{cs0} to obtain
\begin{equation}\label{CSY}
4(H^2+1)(M-1)^2
\leq
\left[\frac{(1-\nu^2)|\nabla\nu|^2}{(H^2+\nu^2)^2}\right].
\end{equation}
Let us estimate this quantity. Recalling that $L\nu=0$, integration by parts gives
\begin{equation}\label{inteparti1}
\int_\Sigma
\frac{(1-\nu^2)|\nabla\nu|^2}{(H^2+\nu^2)^2}\,dA
=-\int_\Sigma q_H(\nu)\Delta\nu\,dA
=\int_\Sigma P\nu \,q_H(\nu)\,dA.
\end{equation}
where
\begin{equation}\label{decu}
q_H(t):=\int_0^t\frac{1-s^2}{(H^2+s^2)^2}\,ds.
\end{equation} 
In order to compute the last integral in \eqref{inteparti1}, we will apply again \eqref{forevenP}, this time to $f(t):=tq_H(t)$, which is smooth and even. Denoting $F(t):=\int_0^t f(s)ds$, as usual, we will use the identities
\[
t(H^2+t^2)q_H(t)-2t F(t)
=(H^2+1)(f_H(t)-1)-t^2
\]
and
\[
F(1)
=\frac12\int_0^1\frac{(1-t^2)^2}{(H^2+t^2)^2}\,dt.
\]
Both identities follow by integration by parts, using \eqref{decu}. Thus, \eqref{inteparti1} and  \eqref{forevenP} 
give 
\begin{equation}\label{intefi1}
\begin{aligned}\int_\Sigma
\frac{(1-\nu^2)|\nabla\nu|^2}{(H^2+\nu^2)^2}\,dA
=&\ 4\int_\Sigma
\bigl((H^2+1)(f_H(\nu)-1)-\nu^2\bigr)\,dA\\
{} & -2\pi\chi(\Sigma)
\int_0^1\frac{(1-t^2)^2}{(H^2+t^2)^2}\,dt.
\end{aligned}
\end{equation}
We finally bound these two terms. For the first one, we note the elementary inequality
\begin{equation}\label{fiti5}
H^2(f_H(t)-1)=Ht\arctan\frac tH\leq t^2, \hspace{0.5cm} \forall t\in [-1,1],
\end{equation}
which implies
\begin{equation}\label{fiti2}
(H^2+1)(f_H(t)-1)-t^2\leq f_H(t)-1.
\end{equation} For the second one,
\begin{equation}\label{fiti3}
\int_0^1\frac{(1-t^2)^2}{(H^2+t^2)^2}\,dt
\leq\frac1{H^2}\int_0^1\frac{dt}{H^2+t^2}
=\frac1{H^3}\arctan\frac1H.
\end{equation}
Plugging \eqref{fiti2}, \eqref{fiti3} and $\chi(\Sigma)\leq 0$ into \eqref{intefi1} gives
\begin{equation}\label{fiti1}
\left[\frac{(1-\nu^2)|\nabla\nu|^2}{(H^2+\nu^2)^2}\right]
\leq
4(M-1)-\frac{2\pi\chi(\Sigma)}{H^3|\Sigma|}
\arctan\frac1H.
\end{equation}
Combining \eqref{fiti1} with \eqref{CSY}, we obtain \eqref{meanb2}.
\end{proof}
Once here, we can complete the proof. First, averaging in \eqref{fiti5} and using \eqref{coarsebounds}, we get $H^2(M-1)\leq[\nu^2]<1/2$, that is,
\begin{equation}\label{boundii}
1\leq M < 1+\frac{1}{2H^2}.
\end{equation}
On the other hand, the quadratic polynomial $p_H(x)$ in \eqref{polipe} satisfies
\begin{equation}\label{directcomparison}
p_H(x)<
\frac45\bigl((H^2+1)(x-1)^2-(x-1)\bigr), \hspace{0.5cm} \forall x\in \left[1, 1+\frac{1}{2H^2}\right].
\end{equation}
Hence, inserting \eqref{boundii}-\eqref{directcomparison} and \eqref{meanb2} into \eqref{roughenergy},  we obtain
\[
\begin{aligned}
\frac{Q(v_H,v_H)}{4H^2|\Sigma|}
&\leq p_H(M)-\frac{\pi\chi(\Sigma)}{H^2|\Sigma|}C_H\\
&<
\frac45\bigl((H^2+1)(M-1)^2-(M-1)\bigr)
-\frac{\pi\chi(\Sigma)}{H^2|\Sigma|}C_H\\
&\leq
-\frac{\pi\chi(\Sigma)}{|\Sigma|}
\left(
\frac{2}{5H^3}\arctan\frac1H+\frac{C_H}{H^2}
\right).
\end{aligned}
\]
Recalling the value of the constant $C_H$ in \eqref{CHval}, we have then
$$
\frac{Q(v_H,v_H)}{4H^2|\Sigma|}<\frac{3\pi\chi(\Sigma)}{4H^2|\Sigma|}
\leq0,$$
which gives \eqref{inegen}. This contradicts weak stability because $[v_H]=0$.

Therefore, any closed weakly stable CMC surface immersed in $\Nil$ has genus $\gamma=0$. By the Abresch--Rosenberg theorem
\cite[Theorem 6]{AR05}, it is an embedded rotational sphere.
This proves Theorem \ref{thm:main}. 

%
%

\subsection*{Statement on AI use} The genus reduction strategy used to prove Theorem B via Abresch-Rosenberg and the Torralbo-Urbano test function $f_H(\nu)$ is due to the authors. LLMs were used to find the final test function $v_H$ and complete the technical computations. The resulting arguments were substantially modified and completely rewritten by the authors, who take responsibility for the correctedness of the paper.

\end{document}